\documentclass[12pt,a4paper]{amsart}
\usepackage{graphicx,amsmath,amssymb,amsthm,color} 

\theoremstyle{plain}
\newtheorem{thm}{Theorem}[section]
\newtheorem{theorem}{Theorem}[section]
\newtheorem{lemma}[theorem]{Lemma}

\newtheorem{conjecture}[theorem]{Conjecture}
\newtheorem{proposition}[theorem]{Proposition}

\theoremstyle{definition}

\newcommand\cB{{\mathcal B}}
\newcommand\cA{{\mathcal A}}
\newcommand\cC{{\mathcal C}}
\newcommand\cF{{\mathcal F}}
\newcommand\cG{{\mathcal G}}

\newcommand{\sats}{{\rm sat}^*}
\newcommand{\sat}{{\rm sat}}

\title{Poset saturation of unions of chains}
\author{Shengjin Ji}
\address{School of Mathematics and Statistics, Shandong University of Technology, Zibo, Shandong 255000, China.}
\thanks{ Shengjin  Ji is supported by  Natural Science Foundation of Shandong Province, China (No. ZR2022MA077) and  Postgraduate Education Reform Project of Shandong Province, China (No. SDYKC2023107)}
\email{jishengjin@sdut.edu.cn}

\author{Bal\'azs Patk\'os}
\address{HUN-REN Alfr\'ed R\'enyi Institute of Mathematics} 
\email{patkos@renyi.hu}

\author{Erfei Yue}
\address{HUN-REN Alfr\'ed R\'enyi Institute of Mathematics \& E\"otv\"os University, Budapest, Hungary} 
\thanks{Erfei Yue is supported by ERC Advanced grant GeoScape, No. 882971. }
\email{yef9262@mail.bnu.edu.cn}

\begin{document}

\begin{abstract}
    A family $\cG$ of sets is a(n induced) copy of a poset $P=(P,\leqslant)$ if there exists a bijection $b:P\rightarrow \cG$ such that $p\leqslant q$ holds if and only if $b(p)\subseteq b(q)$. The induced saturation number $\sats(n,P)$ is the minimum size of a family $\cF\subseteq 2^{[n]}$ that does not contain any copy of $P$, but for any $G\in 2^{[n]}\setminus \cF$, the family $\cF\cup \{G\}$ contains a copy of $P$. We consider $\sats(n,P)$ for posets $P$ that are formed by pairwise incomparable chains, i.e. $P=\bigoplus_{j=1}^mC_{i_j}$. We make the following two conjectures: (i) $\sats(n,P)=O(n)$ for all such posets and (ii) $\sats(n,P)=O(1)$ if not all chains are of the same size. (The second conjecture is known to hold if there is a unique longest among the chains.) We verify these conjectures in some special cases: we prove (i) if all chains are of the same length, we prove (ii) in the first unknown general case: for posets $2C_k+C_1$. Finally, we give an infinite number of examples showing that (ii) is not a necessary condition for $\sats(n,P)=O(1)$ among posets $P=\bigoplus_{j=1}^mC_{i_j}$: we prove $\sats(n,(\binom{2t}{t}+1)C_2)=O(1)$ for all $t\ge 1$.
\end{abstract}


\maketitle
\section{Introduction}

We use standard notations: $[n]$ denotes the set $\{1,2\dots,n\}$ of the first $n$ positive integers and $[i,j]$ denotes the interval $\{i,i+1,\dots,j\}$ for $i<j$. For any set $X$, we write $2^X=\{Y:Y\subseteq X\}$ and $\binom{X}{k}=\{Y\subseteq X:|Y|=k\}$. Whenever the underlying set $[n]$ is clear from context, for a set $F\subseteq [n]$, $F^c$ denotes the complement $[n]\setminus F$. For a family $\cF\subseteq 2^{[n]}$, we write $\cF^c$ to denote the family $\{F^c:F\in \cF\}$ of complements.

In extremal set theory, containment patterns are described by the following definition of Katona and Tarj\'an \cite{KatTar83}: a set system $\cG$ is a \textit{not necessarily induced copy} of a poset $P$ (with some abuse of notation we write $P$ for the poset $(P,\leqslant)$ as the ordering will always be clear from context) if there exists a bijection $b:P\rightarrow \cG$ such that $p\leqslant_P q$ implies $b(p)\subseteq b(q)$. If the bijection $b$ satisfies $p\leqslant_P q$ if and only if $b(p)\subseteq b(q)$, then $\cG$ is an \textit{induced copy} of $P$. If a set system $\cF$ does not contain not necessarily induced copies of $P$, then we say $\cF$ is \textit{$P$-free}, and if $\cF$ does not contain induced copies of $P$, then $\cF$ is \textit{induced $P$-free}. The Tur\'an type problem of finding the largest possible size of a(n induced) $P$-free family $\cF\subseteq 2^{[n]}$ has attracted a lot of attention. We refer the interested reader to \cite{GL} and Chapter 7 of the book \cite{GerPat18}. 

We will be interested in the saturation version of the problem: $\cF \subseteq 2^{[n]}$ is \textit{$P$-saturated} if it is $P$-free, but for any $G\in 2^{[n]}\setminus \cF$, the family $\cF\cup \{G\}$ contains a not necessarily induced copy of $P$. The investigation of $\sat(n,P)$, the minimum size
of a $P$-saturated family $\cF\subseteq 2^{[n]}$ was initiated in \cite{Getal} for the chain poset $C_k$ (bounds later improved in \cite{MNS} and then in \cite{MV}). Later, Keszegh et al \cite{KLMPP} showed that $\sat(n,P)\le C_P$ holds for any poset $P$ for a constant $C_P$ that depends only on $P$ and not on the size $n$ of the underlying set. Therefore, most attention is given to \textit{induced $P$-saturated} families, i.e. induced $P$-free families $\cF$ with the property that for any $G\in 2^{[n]}\setminus \cF$ the family $\cF\cup \{G\}$ contains an induced  copy of $P$. The first results on the \textit{induced saturation number} $\sats(n,P)$, the minimum size of an induced $P$-saturated family $\cF\subseteq 2^{[n]}$, were obtained in \cite{Fetal}. Then results for specific posets \cite{I,I2,MSW,MV2} and for classes of posets, notably antichains \cite{Betal,DI} and complete bipartite posets \cite{Liu} have been established, but our understanding of the behavior of $\sats(n,P)$ is far from satisfying. It was observed in \cite{KLMPP} that for any poset $P$ either $\sats(n,P)\le C_P$ or $\sats(n,P)\ge \log_2 n$ holds. This was improved in \cite{FPST} to  $\sats(n,P)\le C_P$ or $\sats(n,P)\ge 2\sqrt{n}$, but as we only know examples of posets with constant or linear saturation number, we do not know whether one can still improve on this dichotomy result. A general polynomial upper bound of $\sats(n,P)\le Cn^{\frac{|P|^2}{2}}$ was shown in \cite{Betal+}, and recently a poset-operation preserving unbounded saturation number was introduced in \cite{IJ}.

At the moment, it is unknown whether there exists a poset $P$ with superlinear induced saturation number. The aim of the present paper is to formulate a conjecture that a certain class of posets will not provide such an example, and to prove this conjecture in some special cases. For two posets $P,P'$, we denote by $P+P'$ the poset obtained by considering the disjoint union of incomparable copies of $P$ and $P'$, and for a positive integer $k$ and poset $P$, we write $kP$ to denote the union of $k$ pairwise disjoint and uncomparable copies of $P$. We shall consider disjoint unions of chains and we propose the following conjecture.

\begin{conjecture}\label{general}
    For any $i_1\ge i_2 \ge \dots\ge i_m$, there exists $C=C(i_1,i_2, \dots,i_m)$ such that $\sats(n,\bigoplus_{j=1}^mC_{i_j})\le C\cdot n$.
\end{conjecture}

If true, Conjecture \ref{general} is best possible, as, by a result of \cite{KLMPP}, $\sats(n,2C_2)=\Theta(n)$. We propose another conjecture that states that if one considers the disjoint union of chains not all of which are of the same length, then the induced saturation number is constant.

\begin{conjecture}\label{constant}
     For any $i_1\ge i_2 \ge \dots\ge i_m$ with $i_1>i_m$, there exists $C=C(i_1,i_2,\dots, i_m)$ such that $\sats(n,\bigoplus_{j=1}^mC_{i_j})\le C$.
\end{conjecture}

Let us have two remarks: on the one hand Conjecture \ref{constant} cannot give a characterization of posets of disjoint union of chains with bounded induced saturation number as in \cite{KLMPP} $\sats(n,3C_2)$, $\sats(n,5C_2)$, and $\sats(n,7C_2)$ were all proved to be bounded by a constant. On the other hand, an important special case of Conjecture \ref{constant} follows from Theorem 3.8 of \cite{KLMPP}.

\begin{thm}[special case of Theorem 3.8 in \cite{KLMPP}]\label{old}
     For any $i_1\ge i_2 \ge \dots\ge i_m$ with $i_1>i_2$, there exists $C=C(i_1,i_2,\dots, i_m)$ such that $\sats(n,\bigoplus_{j=1}^mC_{i_j})\le C$.
\end{thm}

Now we list our results. Our main theorem states that if all chains have the same length, then conjecture \ref{general} holds true.

\begin{thm}\label{equallength}
    For any pair $m$ and $k$ of positive integers, we have $\sats(n,mC_k)\le C\cdot n$ for some $C=C(m,k)$.
\end{thm}

We also prove the following special case of Conjecture \ref{constant}. Because of Theorem \ref{old}, this is the simplest but somewhat general case of Conjecture \ref{constant} that had not been proved before.

\begin{thm}\label{+C1}\
    
        $\sats(n,2C_k+C_1)\le C(k)$ for some constant depending on $k$.
        
\end{thm}

Our final contribution gives an infinite sequence of values of $m$ with the property that $\sats(n,mC_2)$ is bounded by a constant. All these values $m$ are odd and as $\sats(n,2C_2)=\Theta(n)$, one still wonders whether it is true that $\sats(n,mC_2)$ is a constant if and only if $m$ is odd.

\begin{thm}\label{C2}
    For any $t\ge 2$, we have $\sats(n,(\binom{2t}{t}+1)C_2)\le C(t)$ for some constant depending on $t$.
\end{thm}

Note that Theorem \ref{C2} is also true for $t=1$, but that was already proved in \cite{KLMPP} using a different construction.

\medskip

All our results are upper bounds on saturation numbers, so all our results require constructions. Let us have three general remarks.

\begin{itemize}
    \item 
    As all results are about induced saturation number, we will omit the adjective induced and just write copies of a poset $P$ for induced copies.
    \item 
    We will use the fact that if $\cF$ is induced $P$-free and the family $\cG\subseteq 2^{[n]}$ of sets $G$ with $\cF\cup \{G\}$ not containing a copy of $P$, then $\sats(n,P)\le |\cF|+|\cG|$ as one can greedily select an induced $P$-saturated family $\cF'$ containing $\cF$ and contained in $\cF\cup \cG$. So in most of our proofs we will give a family $\cF$ and prove its $P$-free property and show copies in $\cF\cup \{G\}$ with a small number of exceptions (linear for Theorem \ref{equallength} and constant for Theorem \ref{+C1} and Theorem \ref{C2}).
    \item 
    All posets $P$ considered are unions of pairwise incomparable chains. When showing a copy of $P$ in $\cF\cup \{G\}$, we need to define the corresponding chains $\cC_1,\cC_2,\dots,\cC_m$ and for any $i\neq j$ find an \textit{$(i,j)$-witness} $x_{ij}\in B_i\setminus T_j$, where $B_i$ is the bottom set of $\cC_i$ and $T_j$ is the top set of $\cC_j$.
\end{itemize}

\section{Constructions, proofs}

Throughout this section, we will assume that $n$ is large enough with respect to $m$ and $k$.

\subsection{Construction for $mC_k$} In this subsection we define a construction that shows Theorem \ref{equallength}. 

\medskip

For any~$s\in[m-1]$, let us define the chain
\[
\mathcal{C}_s:=\{\{s\}\cup [m,t]: t\in [m+k-1,n]\}\cup\{\{s\}\}.
\]
For any~$z\in[2,k]$, let
\[
\mathcal{F}_z:=\left\{F\in\binom{[n]}{z}: |F\cap [m+k-3]|\geq z-1\right\}.
\]
For any~$z\in[0,k-2]$, let
\[
\mathcal{F}^z:=\left\{[n]\setminus F: F\in\binom{[m+k-3]}{m+k-3-z}\right\}.
\]
Finally, let
\[
\mathcal{F}=\cF(m,k):=\left(\bigcup_{s=1}^{m-1}\mathcal{C}_s\right)\cup\left(\bigcup_{z=2}^k\mathcal{F}_z\right)\cup\left(\bigcup_{z=0}^{k-2}\mathcal{F}^z\right)
\cup\{\emptyset,[n]\}.
\]

First we show that there is no copy of ~$mC_k$ in~$\mathcal{F}$. Clearly, no copy $\cG$ of $mC_k$ can contain $\emptyset$ or $[n]$ as they are comparable to all other sets. Furthermore, observe that, as every $\cC_s$ forms a chain, if ~$G\in \cG\cap\mathcal{C}_s$,
then sets of ~$\mathcal{C}_s$ cannot belong to chains of $\cG$ other than that of $G$, so we can construct at most~$m-1$ ~$C_k$'s that use elements from~$\mathcal{C}_s$'s. 
On the other hand, we claim that $\bigcup_{z=2}^k\cF_z\cup \bigcup_{z=0}^{k-2}\cF^z$ is $C_k$-free. Note that for any $F\in \cF_z$, we have $|F\cap [m+k-3]|\in \{z-1,z\}$, while for any $F\in \cF^z$, we have $|F\cap [m+k-3]|=z$. With this in mind, we can observe that
\begin{enumerate}
    \item 
    $A\in \bigcup_{z=2}^k\cF_z$ cannot contain any $B\in \bigcup_{z=0}^{k-2}\cF^z$,
    \item 
    if $A\in \cF_z,B\in \cF_{z'}$ with $A\subseteq B$, then $z< z'$ holds,
    \item 
    if $A\in \cF^z,B\in \cF^{z'}$ with $A\subseteq B$, then $z< z'$ holds,
    \item 
    if $A\in \cF_z,B\in \cF^{z'}$ with $A\subseteq B$, then $z\le z'+1$ holds.
\end{enumerate} 
Therefore, by (2) and (3), there is no $C_k$ in $\bigcup_{z=2}^k\cF_z$ or in $\bigcup_{z=0}^{k-2}\cF^z$, and a chain containing sets both from $\bigcup_{z=2}^k\cF_z$ and $\bigcup_{z=0}^{k-2}\cF^z$, by (4), can contain one set from each of $\cF_2,\cF_3,\dots,\cF_j,\cF^{j-1},\cF^j,\dots,\cF^{k-2}$ for some $2\le j\le k-1$, so a total of at most $k-1$ sets.

\medskip

We want to show that $\cF\cup \{G\}$ contains a copy of $mC_k$ for all but $O(n)$ sets $G\in 2^{[n]}\setminus \cF$. To reduce the complexity of the discussion, we introduce the following concept.
Fix arbitrary~$A_0,A_1,\ldots,A_k\subseteq [m+k-3]$ such that~$|A_i|=i$ and~$A_i\subseteq A_{i+1}$.
For any~$a\not\in [m+k-3]$, let~$A_i(a)=A_i\cup\{a\}$. Let~$B_i=[n]\setminus([m+k-3]\setminus A_i)$. 
Then we have the following diagram of triple chains, in which the symbol~$A\rightarrow B$ means~$A\subseteq B$. 
\[
\begin{array} {ccccccccccccccc}
 & & & & A_2 & \rightarrow & A_3 & \rightarrow & \cdots & \rightarrow & A_{k-2} & \rightarrow & A_{k-1} & \rightarrow & A_k \\ 
 & & & & \downarrow & & \downarrow & & & & \downarrow & & \downarrow & & \\
 & & A_1(a) & \rightarrow & A_2(a) & \rightarrow & A_3(a) & \rightarrow & \cdots & \rightarrow & A_{k-2}(a) & \rightarrow & A_{k-1}(a) & & \\
 & & \downarrow & & \downarrow & & \downarrow & & & & \downarrow & & & &  \\
 B_0 & \rightarrow & B_1 & \rightarrow & B_2 & \rightarrow & B_3 & \rightarrow & \cdots & \rightarrow & B_{k-2} & & & & \\
\end{array}
\]

Observe that all sets in the above diagram belong to $\cF$: the $A_i$'s and $A_i(a)$'s to $\cF_i$ and $\cF_{i+1}$, while the $B_i$'s to $\cF^i$. Note also that there are several~$C_{k-1}$ in this diagram, but no~$C_k$. 

\medskip

For any~$G\in 2^{[n]}\setminus\mathcal{F}$, let us write ~$|G|=i$,~$|G\cap[m-1]|=j$, and~$|G\cap[m+k-3]|=l\geq j$. 
Our aim is to include~$G$ somewhere into the diagram above, to introduce a new~$C_k$ and find $m-1$ further $C_k$'s. We proceed by a case analysis. First we assume $k\ge 3$.

\medskip

\textbf{Case 1.}~$k+1\leq i\leq n-m+2$. 

\textbf{Case 1.1.}~$j\geq 1$. 

In this case, either $G=\{s\}\cup [m,m+i-2]$, or there exists an~$a\in[m,m+i-2]\setminus G$. The first type of sets belong to $\cup_{s=1}^{m-1}\cC_s \subseteq \cF$, so we can assume the existence of $a$. For any~$s\in[m-1]$, we begin with the set~$\{s,a\}$. 
We consider a maximal chain $H_1\subseteq H_2\subseteq \dots \subseteq H_{|[m,m+k-3]\setminus\{a\}|}$ of~$[m,m+k-3]\setminus\{a\}$ and obtain $\{s,a\}\subseteq \{s,a\} \cup H_1 \subseteq \dots \{s,a\}\cup H_{|[m,m+k-3]\setminus\{a\}|}$. Here we used $k\ge 3$ as otherwise there is no set $[m,m+k-3]$. If~$a\not\in[m,m+k-3]$, this chain is of length~$k-1$. All $H_i$ belong to some $\cF_z$.
If~$a\in[m,m+k-3]$, this chain is of length~$k-2$ with its largest set~$\{s\}\cup[m,m+k-3]$. We can add~$\{s\}\cup[m,m+k-2]\in \cF_k$ to it to obtain a chain of length~$k-1$. 
Finally, for both cases, we add~$\{s\}\cup[m,m+i-2]\in \cC_s$ to obtain a chain of length~$k$. 
We can do this because we have~$i\geq k+1$. These chains are pairwise incomparable, the $(s,s')$-witness is $s$. We obtained a copy of ~$(m-1)C_k$ and we need one more ~$C_k$. 

\smallskip

\textbf{Case 1.1.1.}~$l\geq k-1$. 

If~$j\geq 2$, let~$A_2$ be a set of two elements taken from~$G\cap [m-1]$.
We choose a set $H$ of ~$k-3$ arbitrary elements in~$(G\cap[m+k-3])\setminus A_2$ (it exists as $l\ge k-1$ and $k\ge 3$), 
and consider a maximal chain from $A_2$ to $A_2\cup H$: ~$A_2\subseteq A_3\subseteq\cdots\subseteq A_{k-1}=A_2\cup H$. 
Fixing an arbitrary~$b\in G\setminus A_{k-1}$ (it exists as $i\ge k+1$), we get a chain of length~$k$ as follows: 
\[
A_2\subseteq A_3\subseteq\cdots\subseteq A_{k-1}\subseteq A_{k-1}(b)\subseteq G. 
\]
This chain is incomparable to the previous $m-1$ copies of $C_k$ as all its sets contain at least two elements from $[m-1]$ and none of them contains $a$.

If~$j=1$, then there is a~$b\in G\setminus[m,m+i-2]$ as sets of the form $\{s\}\cup [m+i-2]$ are in $\cC_s\subset \cF$. Let~$A_1=G\cap[m-1]$, we choose a set $H$ of ~$k-2$ arbitrary elements in~$(G\cap[m+k-3])\setminus A_1$, 
and consider a maximal chain from~$A_1$ to $A_1\cup H$: ~$A_1\subseteq A_2\subseteq\cdots\subseteq A_{k-1}=A_1\cup H$. Then we get a chain of length~$k$ as follows: 
\[
A_1(b)\subseteq A_2(b)\subseteq\cdots\subseteq A_{k-1}(b)\subseteq G.
\]
 This chain is incomparable to all previous chains as all sets contain $b$ and none of them contains $a$, so the $(m,s)$-witness is $b$, the $(s,m)$-witness is $a$.

\smallskip

\textbf{Case 1.1.2.}~$l\leq k-2$. 

Let~$A_l=G\cap[m+k-3]$. For any~$b\in G\setminus A_l$, we have~$A_l\subseteq A_l(b)\subseteq G \subseteq B_l$. 
We can consider a $(k-2-l)$-subset $H$ of $[m+k-3]\setminus (\{a\}\cup A_l)$ and define $A_{k-2}=A_l\cup H$ to get a chain~$B_l\subseteq\cdots\subseteq B_{k-2}=[n]\setminus ([m+k-3]\setminus A_{k-2})$. 

If~$j\geq 2$, we can fix two elements in~$A_l$ to form $A_2$, and get a maximal chain from ~$A_2$ to $A_l$: ~$A_2\subseteq A_3\subseteq\cdots\subseteq A_l$. Then we have a chain of length~$k$ as follows:
\[
A_2\subseteq\cdots\subseteq A_l\subseteq A_l(b)\subseteq G\subseteq B_l\subseteq\cdots\subseteq B_{k-2}. 
\]
This is incomparable to all previous chains as $a\notin B_{k-2}$, and all sets of this chain contain at least two elements of $[m-1]$.

If~$j=1$, then there exists a~$b\in G\setminus [m,m+i-2]$, as all sets of the form $\{s\}\cup [m,m+i-2]$ belong to $\cC_s\subseteq \cF$. Let~$A_1=G\cap[m-1]$, 
we consider a maximal chain from $A_1$ to ~$A_l=G\cap[m+k-3]$: ~$A_1\subseteq\cdots\subseteq A_l$. Then we have a chain of length~$k$ as follows:
\[
A_1(b)\subseteq\cdots\subseteq A_l(b)\subseteq G\subseteq B_l\subseteq\cdots\subseteq B_{k-2}
\]
 This chain is incomparable to all previous chains as all sets contain $b$ and none of them contains $a$, so the $(m,s)$-witness is $b$, the $(s,m)$-witness is $a$.

\textbf{Case 1.2.}~$j=0$. 

For any~$s\in[m-1]$, we have the following chain of length at least~$k$ from $\cC_s$: 
\[
\{s\}\subseteq\{s,m\}\subseteq\{s,m,m+1\}\subseteq\cdots\subseteq\{s\}\cup[m,m+i-2]. 
\]
These chains are incomparable as the $(s,s')$-winess is $s$. Now we have~$m-1$ copies of $C_k$, and we need one more. Let us write~$b=\max G$. As $j=0$, we have ~$b\geq m+i-1$.

\textbf{Case 1.2.1.}~$l\geq 1$. 

Let~$A_l=G\cap[m+k-3]$. Fixing any $x\in A_l$ and writing $A_1=\{x\}$, we can consider a maximal  chain~$A_1\subseteq\cdots\subseteq A_l$. 
On the other hand, let~$B_{k-2}=[n]\setminus[m-1]$. We can consider a maximal chain from ~$B_l=[n]\setminus ([m+k-3]\setminus A_l)$ to~$B_{k-2}$. 
Then we have a chain of length~$k$ as follows:
\[
A_1(b)\subseteq\cdots\subseteq A_l(b)\subseteq G\subseteq B_l\subseteq\cdots\subseteq B_{k-2}. 
\]
This is incomparable to all previous chains as $B_{k-2}\cap [m-1]=\emptyset$ and so the $(s,m)$-witness is $s$, while $b\notin \{s\}\cup [m,m+i-2]$ for all $s\in [m-1]$ so $(m,s)$-witness is $b$.

\textbf{Case 1.2.2.}~$l=0$. 

Let~$B_0=[n]\setminus[m+k-3]$, and~$B_{k-2}=[n]\setminus[m-1]$. Then we can consider a maximal chain from ~$B_0$ to~$B_{k-2}$, and obtain a chain of length~$k$ as follows:
\[
G\subseteq B_0\subseteq\cdots\subseteq B_{k-2}. 
\]
This is incomparable to the previous chains again as $B_{k-2}\cap [m-1]=\emptyset$ and so the $(s,m)$-witness is $s$ and $b\notin \{s\}\cup [m,m+i-2]$ for all $s\in [m-1]$ and so the $(m,s)$-witness is $b$.

\medskip

\textbf{Case 2.}~$i\leq k$. 

\textbf{Case 2.1.}~$\max G\leq m+k-1$. 

The number of these sets is less than $2^{m+k-1}=O(1)$, so these sets can be considered as exceptional sets. 

\textbf{Case 2.2.}~$\max G\geq m+k$. 


\textbf{Case 2.2.1.}~$j=0$. 

We obtain the first~$m-1$ chains as we did in Case 1.2: 
for any~$s\in [m-1]$, we have the following chain of length ~$k$: 
\[
\{s\}\subseteq\{s,m\}\subseteq\{s,m,m+1\}\subseteq\cdots\subseteq\{s\}\cup[m,m+k-2]. 
\]
These are pairwise incomparable as the $(s,s')$-witness is $s$.

If~$l\geq 1$, let~$b=\max G$. We obtain the last chain similar to Case 1.2.1. With $A_l:=G\cap [m+k-3]$ and $A_1=\{x\}$ for some $x\in A_l$, we consider a maximal chain $A_1\subseteq A_2 \subseteq \dots \subseteq A_l$ and for $B_l=[n]\setminus ([m+k-3]\setminus A_l)$ and $B_{k-2}=[n]\setminus [m-1]$ a maximal chain from $B_l$ to $B_{k-2}$, and then define
\[
A_1(b)\subseteq\cdots\subseteq A_l(b)\subseteq G\subseteq B_l\subseteq\cdots\subseteq B_{k-2}.
\]
This is incomparable to the first $m-1$ chains, because $b\ge m+k$ so the $(m,s)$-witness is $b$, while $B_{k-2}\cap [m-1]=\emptyset$ and so the $(s,m)$-witness is $s$.

If~$l=0$, we have the last chain as in Case 1.2.2. With $B_0=[n]\setminus [m+k-3]$ and $B_{k-2}=[n]\setminus [m-1]$, we define 
\[
G\subseteq B_0\subseteq\cdots\subseteq B_{k-2}.
\]
The $(s,m)$-witness is $s$, while the $(m,s)$-witness is $\max G\ge m+k$.

\textbf{Case 2.2.2.}~$j\geq 1$. 

Then, as $i\le k$, there exists an~$a\in[m,m+k-2]\setminus G$. For any~$s\in[m-1]$, we have a chain: 
\[
\{s,a\}\subseteq\cdots\subseteq\{s\}\cup[m,m+k-2]\subseteq\{s\}\cup[m,m+k-1]. 
\]
These are pairwise incomparable as the $(s,s')$-witness is $s$. The construction of the last~$C_k$ that contains~$G$ is as in Case 2.2.1 when $l\ge 1$: we pick $A_1=\{x\}$ for some $x\in A_l:=G\cap [m+k-3]$ and a maximal chain $A_1\subseteq A_2 \subseteq \dots \subseteq A_l\subseteq \dots \subseteq A_{k-2}$ with $A_{k-2}\subseteq [m+k-3]\setminus \{a\}$ and consider the chain 
\[
A_1(b)\subseteq\cdots\subseteq A_l(b)\subseteq G\subseteq B_l\subseteq\cdots\subseteq B_{k-2}
\]
with $b$ still being $\max G\ge m+k$. The $(m,s)$-witness is $b$, while the $(s,m)$-witness is $a$.
\medskip

\textbf{Case 3.}~$i\geq n-m+3$. 

In this case,~$G$ misses at most~$m-3$ elements of $[n]$, so~$j\geq 2$, ~$l\geq k$. 
Let~$A_2$ be a set of two elements in~$G\cap [m-1]$, and $A_k\supseteq A_2$ be a $k$-subset of $G\cap [m+k-3]$. We consider a maximal chain from $A_2$ to $A_k$ and obtain a chain of length~$k$ as follows: 
\[
A_2\subseteq\cdots\subseteq A_k\subseteq G. 
\]
We still need to define the other $m-1$ chains of length $k$. Let~$G'=[n]\setminus G$. 

\textbf{Case 3.1.}~$a=\min (G'\setminus [m-1])\leq n-k+1$. 

For any~$s\in[m-1]$, we have a chain of length~$k$: 
\[
\{s\}\cup[m,a]\subseteq\cdots\subseteq\{s\}\cup[m,a+k-1]. 
\]
The $(s,s')$-witness is $s$, the $(s,m)$-witness is $a$, the $(m,s)$-witness is an element of $A_2\setminus \{s\}$.

\textbf{Case 3.2.}~$G'\subseteq[m-1]$ or~$\min (G'\setminus[m-1])\geq n-k+2$. 

The number of such sets is at most $2^{m+k-2}$, so we can add these to the exceptions.

\medskip

We still need to consider the $k=2$ case. The construction is
\[\mathcal{F}=\left(\cup_{s=1}^{m-1}\right)\mathcal{C}_s\cup\mathcal{F}_2\cup\mathcal{F}^0\cup\{\emptyset,[n]\}
,\]
$\mathcal{F}_2=\binom{[n]}{2}\setminus \binom{[m,n]}{2},\mathcal{F}^0=\{[m,n]\}$,
$
\mathcal{C}_s=\{\{s\}\cup[m,t]: t\in[m+1,n]\}\cup\{\{s\}\}.
$

\smallskip

For any~$G\in 2^{[n]}\setminus\mathcal{F}$, let~$|G|=i$ and~$|G\cap[m-1]|=j$. As there are only finitely many sets containing $[m,n]$ and there are only finitely many sets contained in $[m]$, we can add these to the exceptional sets and assume there exists $a\in [m,n]\setminus G$.

\smallskip


\textbf{Case 1.} $j\geq 2.$

Then for every~$s\in [m-1]$, we have a chain~$\{s,a\}\subseteq\{s\}\cup[m,n]$, 
so~$m-1$ chains altogether. For the last chain, as $j\ge 2$, there exists a 2-subset $H$ of $G\cap [m-1]$ and~$H\subseteq G$ is the last chain. The $(s,s')$-witness is $s$, the $(s,m)$-witness is $a$, $H$ contains two elements from $[m-1]$ while sets in the first $m-1$chains intersect $[m-1]$ in one element.

\textbf{Case 2.} $j=1.$

Note that~$|G\cap[m,n]|=i-1=|[m,m+i-2]|$. As $\{s\}\cup [m,m+i-2]\in \mathcal{C}_s\subseteq \mathcal{F}$ and $G\notin \mathcal{F}$, there exist~$x\in[m,m+i-2]\setminus G$ and ~$y\in G\cap[m+i-1,n]$.
Then for every~$s\in[m-1]$, we have a chain~$\{s,x\}\subseteq\{s\}\cup[m,m+i-2]$, and the last chain is~$(G\cap[m-1])\cup\{y\}\subseteq G$. The $(s,s')$-witnesses is $s$, the $(s,m)$-witness is $x$, the $(m,s)$-witness is $y$.

\textbf{Case 3.} $j=0.$

For every~$s\in [m-1]$, we have a chain~$\{s\}\subseteq\{s,m\}$, and we can simply take the last chain by~$G\subseteq [m,n]$. By our assumption $G\not\subseteq [m]$, there exists $x\in G\setminus [m]$. The $(s,s')$, $(s,m)$, and $(m,s)$-witnesses are $s$, $s$, and $x$, respectively.

\subsection{A constant size construction for $mC_2$ with $m=\binom{2t}{t}+1$}
In Section 2.1, we established a linear upper bound on~$\sats(n,mC_k)$ for all $m$ and $k$. 
In \cite{KLMPP}, constructions yielding constant upper bounds for~$\sats(n,3C_2),\sats(n,5C_2)$ and~$\sats(n,7C_2)$ are presented.
By proving Theorem \ref{C2}, we will show that there is an infinite number of values of $m$ for which $\sats(n,mC_2)$ is bounded by a constant depending only on $m$.
Let~$m=\binom{2t}{t}+1$, for a positive~$t$.
Our construction is as follows. For any ground set~$[n]$, let
\[
\mathcal{F}_0=\{\emptyset\}\cup\binom{[2t+1]}{t}\cup\left\{F\in\binom{[2t+1]}{t+1}: 2t+1\in F\right\}\cup\binom{[2t]}{t+2},
\]
and
\[
\mathcal{F}_1=\{F^c: F\in\mathcal{F}_0\}.
\]
Then we claim that ~$\mathcal{F}=\mathcal{F}_0\cup\mathcal{F}_1$ is $mC_2$-saturated.

\medskip

First of all, we need to prove that~$\mathcal{F}$ is~$mC_2$-free. To prove that, we need the following concept and classical result.
Suppose~$\mathcal{P}=\{(X_i,Y_i): i\in[m]\}$ is a family of pairs of sets, where~$X_i\cap Y_i=\emptyset$.
Then~$\mathcal{P}$ is called a~\emph{Bollob\'as~system} if~$X_i\cap Y_j\neq\emptyset$, for all $1\le i\neq j\le m$,
and a \emph{skew~Bollob\'as~system} if~$X_i\cap Y_j\neq\emptyset$, for all $1\le i<j\le m$.
Note that a~Bollob\'as~system is always a skew~Bollob\'as~system. 
The following lemma was first proved by~Bollob\'as \cite{Sets} only for Bollobás systems, and later strengthened by Lov\'asz and by Frankl to skew Bollob\'as systems.

\begin{thm}[Frankl~\cite{Skew}, Lov\'asz~\cite{L}]\label{Th:skew-uniform}
Let~$\mathcal{P}=\{(X_i,Y_i): i\in[m]\}$ be a skew~Bollob\'as~system, where~$|X_i|\leq a,|Y_i|\leq b$ for any~$i$.
Then~$m\leqslant\binom{a+b}{a}$.
\end{thm}

\begin{proof}[Proof of the fact that ~$\cF$ is~$mC_2$-free.]
Let~$\mathcal{C}=\{B_1,\ldots,B_k,T_1,\ldots,T_k\}$ be a copy of~$kC_2$ in~$\mathcal{F}$, where~$B_i\subseteq T_i$ for each~$i$. 
We need to show that~$k\leq\binom{2t}{t}$. Note that~$B_i\subseteq T_i$ means~$B_i\cap T_i^c=\emptyset$, 
and the condition that $\cC$ is an induced copy of $kC_2$ implies that~$B_i\not\subseteq T_j$, and so~$B_i\cap T_j^c\neq\emptyset$ whenever~$i\neq j$. 
Hence~$\{(B_i,T_i^c): i\in[k]\}$ forms a~Bollob\'as~system. 
Our target is for each~$i$, mapping~$B_i$ to some~$X_i$, and~$T_i$ to some~$Y_i$, such that~$\{(X_i,Y_i): i\in[k]\}$ is still a skew~Bollob\'as~system, 
and~$|X_i|\leq t, |Y_i|\leq t$ for every~$i\in[k]$. Note that each~$B_i$ must belong to one of the following families: 
\begin{multline*}
\mathcal{B}_1=\binom{[2t]}{t}, \mathcal{B}_2=\left\{F\in\binom{[2t+1]}{t}: 2t+1\in F\right\}, 
\mathcal{B}_3=\left\{F: 2t+1\not\in F, F^c\in\binom{[2t+1]}{t+1}\right\}, \\
\mathcal{B}_4=\left\{F\in\binom{[2t+1]}{t+1}: 2t+1\in F\right\}, \mathcal{B}_5=\left\{F: F^c\in\binom{[2t]}{t+2}\right\}. 
\end{multline*}
We may adjust the ordering of the~$C_2$'s such that for any~$i<j$, if~$B_i\in\mathcal{B}_u$ and~$B_j\in\mathcal{B}_v$, then~$u\leq v$. 
If~$B_i$ is in~$\mathcal{B}_1$ or~$\mathcal{B}_2$, let~$X_i=B_i$, and~$Y_i=T_i^c\cap[2t]$. 
Otherwise let~$X_i=B_i\cap[2t]$ and~$Y_i=T_i^c$. 

An ~$X_i$ has size $t$ if $B_i\in \cup_{i=1}^4\cB_i$ and size $t-2$ if $B_i\in \cB_5$. Let us consider the ~$Y_i$'s: if $B_i\in \cB_1$, then as $B_i\cap T_i^c=\emptyset$, we have $|Y_i|\le t$. If $B_i\in \cB_2$, then $T_i\in \cB_4\cup \cB_1^c$, and in both cases $|T_i^c\cap [2t]|=t$. If $B_i\in \cB_3$, then $T_i\in \cB_1^c\cup \cB_2^c$ so $|Y_i|=t$. If $B_i \in \cB_4\cup \cB_5$, then $T_i \in \cB_1^c$ and so $|Y_i|=t$.

\medskip
Finally, for every~$i<j$, we need to check~$X_i\cap Y_j\neq\emptyset$. We proceed by a case analysis.

\medskip

\textbf{Case 1.}~$B_i\in\mathcal{B}_1\cup\mathcal{B}_2$, and~$B_j\in\mathcal{B}_3\cup\mathcal{B}_4\cup\mathcal{B}_5$. 
In this case, we simply have~$X_i\cap Y_j=B_i\cap T_j^c\neq\emptyset$. 

\medskip

\textbf{Case 2.}~$B_i,B_j\in\mathcal{B}_1\cup\mathcal{B}_2$. Then~$X_i=B_i, Y_j=T_j^c\cap[2t]$. 

\smallskip

\textbf{Case 2.1.}~$B_i\in\mathcal{B}_1$. Then~$B_i\subseteq [2t]$. Hence
\[
X_i\cap Y_j=B_i\cap T_j^c\cap [2t]=B_i\cap T_j^c\neq\emptyset. 
\]

\smallskip

\textbf{Case 2.2.}~$B_i\in\mathcal{B}_2$. Then~$B_j\in\mathcal{B}_2$, and~$2t+1\in B_j\subseteq T_j$. Hence~$2t+1\not\in T_j^c$, and
\[
X_i\cap Y_j=B_i\cap T_j^c\cap [2t]=B_i\cap T_j^c\neq\emptyset. 
\]

\medskip

\textbf{Case 3.}~$B_i,B_j\in\mathcal{B}_3\cup\mathcal{B}_4\cup\mathcal{B}_5$. Then~$X_i=B_i\cap[2t], Y_j=T_j^c$. 

\smallskip

\textbf{Case 3.1.}~$B_i,B_j\in\mathcal{B}_3$. Note that~$T_j$ contains all elements no less than~$2t+2$, and~$2t+1\not\in B_i$. 
So~$B_i\not\subseteq T_j$ implies that~$B_i\cap [2t]\not\subseteq T_j$. Hence
\[
X_i\cap Y_j=B_i\cap T_j^c\cap [2t]\neq\emptyset.
\]

\smallskip

\textbf{Case 3.2.}~$B_j\in\mathcal{B}_4\cup\mathcal{B}_5$. In this case,~$T_j$ must lay in $\cB_1^c$. Hence
\[
X_i\cap Y_j=B_i\cap T_j^c\cap [2t]=B_i\cap T_j^c\neq\emptyset.
\]

Therefore~$\{(X_i,Y_i): i\in[k]\}$ is a skew~Bollob\'as~system. 
Then by lemma~\ref{Th:skew-uniform}, we have~$k\leq\binom{2t}{t}$, and so~$\mathcal{F}$ is~$mC_2$-free. 
\end{proof}

\begin{proof}[Proof of~$mC_2$-saturation.]
For every~$G\not\in\mathcal{F}$, we need to prove that there is an induced~$mC_2$ in~$\mathcal{F}\cup\{G\}$. 
As $\cF$ is closed under taking complements, we may assume that~$2t+1\in G$. 
Let~$k=|G\cap[2t]|$. By symmetry, we may assume that~$G\cap[2t]=[k]$. 
Finally, we may assume that there exists an~$x\in[2t+2,n]\setminus G$, because there are only a constant number of exceptions. 

\smallskip

\textbf{Case 1.}~$k\leq t-2$. We have a first $C_2$ as $G\subseteq[t-1,2t]^c$, and for every~$F\in\binom{[2t]}{t}$, we have a chain~$F\subseteq([2t+1]\setminus F)^c$. 
So we have~$m=\binom{2t}{t}+1$ chains altogether. The $(1,j)$-witness is $2t+1$, $F\not\subseteq [t-1,2t]^c$ as $|[2t]\cap [t-1,2t]^c|=t-2$ and the other chains are incomparable as $[2t]\cap ([2t+1]\setminus F)^c)=F$.

\smallskip

\textbf{Case 2.}~$k=t-1$. Then we have first a chain~$[t-1]\cup\{2t+1\}\subseteq G$ and for every~$F\in\binom{[2t]}{t}$, we have a chain~$F\subseteq([2t+1]\setminus F)^c$.
So we have~$m=\binom{2t}{t}+1$ chains altogether. The $(1,j)$-witness is $2t+1$, $F\not\subseteq G$ as $|[2t]\cap  G|=t-1$ and the other chains are incomparable as $[2t]\cap ([2t+1]\setminus F)^c)=F$. 

\smallskip

\textbf{Case 3.}~$t\leq k\leq t+1$. We have a first chain~$B_1=[t]\cup\{2t+1\}\subseteq G$, and a second one~$B_2=[t+1,2t+1]^c\subseteq (\{t+1\}\cup[t+3,2t+1])^c=T_2$. The $(1,2)$-witness is $2t+1$ and the $(2,1)$-witness is $x$.
For every~$F\in\binom{[2t]}{t}$, and~$F\neq[t+1,2t]$, we have a chain~$[2t+1]\setminus F\subseteq F^c$ if~$t+1\in F$, 
or a chain~$(F\cup\{2t+1\})^c\subseteq F^c$ if~$t+1\not\in F$. So we have~$m=\binom{2t}{t}+1$ chains altogether. Incomparability: 
\begin{itemize}
    \item 
    $[2t+1]\setminus F \not\subseteq G$ if $t+1\in F\neq [t+1,2t]$ as $G\cap [2t]\subseteq [t+1]$, but $[t+2,2t]\setminus F\neq \emptyset$ for all such $F$. Also, $[t]\not\subseteq F^c$ (here we use $F\neq [t+1,2t]$), so $B_1\not\subseteq F^c$.
    \item 
    $[2t+1]\setminus F \not\subseteq T_2$ if $t+1\in F\neq [t+1,2t]$ as witnessed by $2t+1$. Also, $[t]\not\subseteq F^c$ and $[t]\subseteq B_2$.
    \item 
    $(F\cup \{2t+1\})^c\not\subseteq G$ if $t+1\notin F$ as witnessed by $x$, and $[t]\not\subseteq F^c$ as $|F^c\cap [2t]|=t$ and $t+1\in F^c$.
    \item 
    $(F\cup \{2t+1\})^c\not\subseteq T_2$ if $t+1\notin F$ as witnessed by $t+1$. Also, $B_2\not\subseteq F^c$ as in the previous point.
    \item 
    The incomparability of the chains for $F$'s is clear.
\end{itemize}

\smallskip

\textbf{Case 4.}~$k\geq t+2$. We have a first chain~$[t+2]\subseteq G$, and for every~$F\in\binom{[2t]}{t}$, we have a chain~$(F\cup\{2t+1\})^c\subseteq F^c$. 
So we have~$m=\binom{2t}{t}+1$ chains altogether. The $(j,1)$-witness is $x$, and $[t+2]\not\subseteq F^c$ as $|F^c\cap [2t]|=t$, the other chains are incomparable as their restrictions to $[2t]$ are distinct $t$-sets $[2t]\setminus F$. 
\end{proof}

\subsection{Construction for $2C_k+C_1$}

\medskip
In this subsection, we give a construction to prove Theorem \ref{2ckc1}. Before defining our construction, let us introduce some notation.
Let $\mathbb{B}_{k,-}$, $\mathbb{B}_{k,-,-}$ be the posets on $2^{[k]}\setminus \{\emptyset\}$ and $2^{[k]}\setminus \{\emptyset,[k]\}$, respectively, both ordered by inclusion. The minimal elements of these posets are called their \textit{atoms}.

Fix a~$k$, and ~$n\ge 2k$. Let $\cA_1=2^{[k]}\setminus \{\emptyset\},\cA_2=2^{[k+1,2k-1]}\setminus \{\emptyset\},\cA_{1,2}=\{[k]\cup A:A\in \cA_2\}$. Define
\[
\mathcal{F}=\cF_k=\{\emptyset,[n]\}\cup (\cA_1 \cup \cA_2 \cup \cA_{1,2}) \cup  (\cA_1 \cup \cA_2 \cup \cA_{1,2})^c. 
\] 
Then~$|\mathcal{F}|=2^{k+2}-4$. 

We now gather some basic properties of $\mathbb{B}_{k,-},\mathbb{B}_{k,-,-}$ and the construction $\cF$. For a family $\cG$ and a poset $P$, we write $\cG\leqslant P$ if the poset on $\cG$ ordered by inclusion  has a copy in $P$ and we write $\cG=P$ if this poset is isomorphic to $P$.

\begin{lemma}\label{bk-} \
    \begin{enumerate}
        \item 
        $\mathbb{B}_{k,-}$ is $(C_k+C_1)$-free and thus $2\mathbb{B}_{k,-}$ is $(2C_k+C_1)$-free.
        \item 
        $\mathbb{B}_{k+1,-,-}$ is $(2C_k+C_1)$-free.
        \item 
        $\cA_1=\mathbb{B}_{k,-}$,  $\cA_2=\mathbb{B}_{k-1,-}$, $\cA_1\cup \cA_1^c= \mathbb{B}_{k+1,-,-}$.
        \item 
        $\cA_2\cup \cA_{1,2}\cup \cA_2^c\cup \cA_{1,2}^c\leqslant \mathbb{B}_{k+1,-,-}$.
        \item 
        $\cA_2\cup \cA_{1,2}^c= \mathbb{B}_{k,-,-}$.
    \end{enumerate}
\end{lemma}

\begin{proof}
    To see (1), observe that a chain of length $k$ in $\mathbb{B}_{k,-}$ must contain $[k]$ which is comparable to all other sets.

    To see (2), observe that $k$-chains in $\mathbb{B}_{k+1,-,-}$ must contain sets of size 1 through $k$, so for two such incomparable chains there must exist $x,y$ with one chain going from $\{x\}$ to $[k+1]\setminus \{y\}$, while the other going from $\{y\}$ to $[k+1]\setminus \{x\}$. But every set is either comparable to $\{x\}$ or to $[k+1]\setminus \{x\}$.

    (3) is trivial.
    
    To see (4), consider $[k],\{k+1\},\{k+2\},\dots,\{2k-1\},[2k,n]$ as atoms, while to see (5), consider $\{k+1\},\{k+2\},\dots,\{2k-1\},[2k,n]$ as atoms.
\end{proof}

The following theorem implies Theorem \ref{+C1}.

\begin{thm}\label{2ckc1}
 $\mathcal{F}$ is $2C_k+C_1$-free, and for all but a constant number of $G\in 2^{[n]}\setminus \cF$ the family $\cF\cup \{G\}$ contains a copy of $2C_k+C_1$.    
\end{thm}
   
\begin{proof} 
We first check that $\cF$ is ~$(2C_k+C_1)$-free. Suppose towards a contradiction that $\cG\subseteq \cF$ is a copy of $2C_k+C_1$. Clearly, $\emptyset,[n]\notin \cG$ as they are comparable to all other sets in $2^{[n]}$. Assume first $(\cA_1\cup \cA_1^c)\cap \cG=\emptyset$. Then $\cG\subseteq \cA_2\cup \cA_{1,2}\cup \cA_2^c\cup \cA_{1,2}^c$ contradicting ~Lemma \ref{bk-} (2) and (4).

Assume next $\cG$ contains a set from only one of $\cA_1$ and $\cA_1^c$. Without loss of generality we may assume that $A\in \cA_1\cap \cG$. Then $A$ is comparable to all sets of $\cA_{1,2}\cup \cA_2^c$. The two chains of $\cG$ not containing $A$ 
cannot be completely in $\cA_2\cup \cA_{1,2}^c$ because of Lemma \ref{bk-} (5). So one of these two chains must contain an $A'$ from $\cA_1$. As $A$ is incomparable to all sets in $\cA_2\cup \cA_{1,2}^c$, we must have that the chain of $\cG$ containing $A$ is completely in $\cA_1$ and thus the chain of $A$ and $A'$ would form a copy of $C_k+C_1$ in $\cA_1$ contradicting ~Lemma \ref{bk-} (1) and (3).

Finally, assume there exist $A\in \cA_1\cap \cG$ and $B\in \cA_1^c\cap \cG$. Every $H \in \cA_2\cup \cA_{1,2} \cup \cA_2^c\cup \cA_{1,2}^c$ is comparable to either $A$ or $B$, so the third chain in $\cG$ (the one that does not contain neither $A$ nor $B$) must lie in $\cA_1\cup \cA_1^c$. By symmetry, we may and will assume that there exists a set $A'\in \cA_1\cap \cG$ incomparable to both $A,B$. As $A'$ is comparable to all sets in $\cA_{1,2}\cup \cA_2^c$ and $A$ is incomparable to all sets in $\cA_2\cup \cA_{1,2}^c$, the chain of $A$ must lie in $\cA_1\cup \cA_1^c$. By Lemma \ref{bk-} (1) and (3), the chain of $A$ must contain a set from $\cA_1^c$, otherwise we would have $C_k+C_1$ in $\cA_1$. So we can assume that $A,B$ belong to the same chain. But then the third chain (that does not contain $A,B,A'$) is incomparable to both $A,B$, so lies in $\cA_1\cup \cA_1^c$ and thus $\cG \subseteq \cA_1\cup \cA_1^c$ contradicting Lemma \ref{bk-} (2) and (3).

\medskip

    Now we are going to prove that for arbitrary~$G\in 2^{[n]}\setminus\mathcal{F}$, the family ~$\mathcal{F}\cup \{G\}$ contains a copy of $2C_k+C_1$. 
Note that~$G$ is incomparable with at least one of~$[k]$ and~$[k]^c$ as $\emptyset, [n]\in \cF$. As $\cF$ is closed under complementation, we may assume that~$G$ is incomparable with~$[k]^c$. 

\smallskip

\textbf{Case 1.} $G$ is incomparable with $[k]$, so there exists $a\in [k]\setminus G$, $b \in G\setminus [k]$.

\smallskip

Because $G$ and $[k]^c$ are incomparable, there exists $a'\in [k]\cap G$. Let us fix $b'\neq b$ with $b'\in [k+1,2k-1]$ (this can be done as $k\ge 3$).
Then one can take the first $C_k$ in $\cA_1$ from $\{a,a'\}$ to $[k]\cup \{b'\}$, and another $C_k$ starting with $([k]\cup\{b'\})^c\in \cA_{1,2}^c$ and continuing in $\cA_{1}^c$ from $[k]^c$ to $\{a,a'\}^c$, while $G$ can serve as the extra $C_1$. The $(1,2)$- and $(1,3)$-witness is $a$, the $(3,1)$-witness is $b$, the $(2,1)$-witness and the $(3,2)$-witness is $a'$. A $(2,3)$-witness does not exist if $G$ is a superset of $([k] \cup \{b'\})^c$, but there is only a finite number of such exceptional sets.

\smallskip

\textbf{Case 2.} $G$ is comparable with $[k]$ and, as $\cA_1\subseteq \cF$, $[k]\subseteq G$.

\smallskip

If $G$ is comparable with $[2k-1]$, then $[2k-1]\subseteq G$ as all subsets of $[2k-1]$ that contain $[k]$ are in $\cA_{1,2}\subseteq \cF$. As $[n]\in \cF$, there exists $x\in [2k,n]\setminus G$. Consider the following chains: one of length $k$ from $[k]^c$ to $\{b\}^c$ for any $b\in [k]$, these sets are in $\cA_1^c$. Another one of length $k$ from $[k+1]$ to $[2k-1]$ (these sets are in $\cA_{1,2}$) and then $G$. Finally, $[k+1,2k-1]^c\in \cA_2^c$ is a chain of length one. All three chains are incomparable as witnessed by $b$, $x$, and $k+1$, so they form a $2C_k+C_1$ in $\cF\cup \{G\}$.

If $G$ is incomparable with $[2k-1]$, then there exists $x\in [k+1,2k-1]\setminus G$. If further $G$ is comparable with $[k+1,2k-1]^c$, then $[k]\subseteq G\subseteq [k+1,2k-1]^c$ as all supersets of $[k+1,2k-1]^c$ are in $\cA_2^c\subseteq \cF$. Then we can consider the following chains: one of length $k$ from $[k]^c$ to $\{b\}^c$ for any $b\in [k]$ (sets in $\cA_1^c$), another one of length $k$ containing $G$ and then a chain from $[k+1,2k-1]^c$ to $\{x\}^c$ in $\cA_2^c$, and finally the single set $[2k-1]$. These three chains are pairwise incomparable as witnessed by $b$ and $x$, so they form  a $2C_k+C_1$ in $\cF\cup \{G\}$.

Finally, if $G$ is incomparable with both $[2k-1]$ and $[k+1,2k-1]^c$. Then there exist $x\in G\cap [k+1,2k-1]$, $x'\in G\cap [2k,n]$, $y\in [k+1,2k-1]\setminus G$ and $y'\in [2k,n]\setminus G$. The three chains are: a chain of length $k$ as $[k]\cup \{x\}\subseteq\dots \subseteq G\cap [2k-1] \subseteq G \subseteq G \cup [k+1,2k-1]^c \subseteq\dots \subseteq\{y\}^c$ (because of $x'$ the sets $G\cap [2k-1]$ and $G$ are distinct, and because of $y'$, the sets $G$ and $G\cup [k+1,2k-1]^c$ are distinct), another chain of length $k$ from $[k]^c$ to $\{1\}^c$ in $\cA_2^c\cup \cA_{1,2}^c$ and $[k]\cup \{y\}$ as a chain of length one. These chains are incomparable as witnessed by 1, $y$, and $x$ so they form a $2C_k+C_1$ in $\cF\cup \{G\}$. 
\end{proof}

\section{Final remarks}

We introduced two conjectures concerning the order of magnitude of $\sats(n,\bigoplus_{j=1}^mC_{i_j})$. While we strongly believe in Conjecture \ref{general}, we would not be that much surprised if a counterexample to Conjecture \ref{constant} were to be found. According to a colleague of ours, who prefers to stay anonymous, a possible such counterexample could be of the form $2C_k+2C_l$.

\smallskip

We tried to generalize the construction from Section 2.3 to the case $P=3C_k+C_1$ (and then to $P=mC_k+C_1$), but all our efforts resulted in families $\cF$ that are either not $P$-free or there are too many sets $G\notin \cF$ that do not yield copies of $P$ in $\cF\cup \{G\}$. However, the case $P=mC_k+C_1$ of Conjecture \ref{constant} is the most disturbing that we could not settle.

\smallskip

Finally, it would be great to define some classes of posets $P$ for which one can determine the order of magnitude of $\sats(n,P)$. One such class could be $\mathbb{B}_n$. The following proposition shows that $\sats(n,\mathbb{B}_3)$ is at most linear in $n$.

\begin{proposition}
    $\sats(n,\mathbb{B}_3)\le 3n-2$.
\end{proposition}

\begin{proof}
    Consider the family $\cF:=\{\emptyset,[n]\}\cup\binom{[n]}{1}\cup \{uv:u\in \{1,2\},v\in [3,n]\}$. That is, $\cF$ contains the empty set, the underlying set, all singletons, and the edge set of the complete bipartite graph $K_{2,n-2}$. Clearly, $|\cF|=1+1+n+2(n-2)=3n-2$. First, we claim that $\cF$ is $\mathbb{B}_3$-free. As there are only four different set sizes in $\cF$, a copy of $\mathbb{B}_3$ should contain $\emptyset,[n]$, three singletons, and 3 pairs that form a triangle. But $K_{2,n-2}$ is triangle-free, so $\cF$ is $\mathbb{B}_3$-free.

    We need to show that for any $G\in 2^{[n]}\setminus \cF$, the family $\cF\cup \{G\}$ contains a copy of $\mathbb{B}_3$. If $|G|=2$, then $G=\{1,2\}$ or $G=\{u,v\}$ with $3\le u,v\le n$. In the former case $2^{[3]}\setminus \{[3]\} \cup \{[n]\}\subseteq \cF\cup \{G\}$, in the latter case $2^{\{1,u,v\}}\setminus \{\{1,u,v\}\}\cup \{[n]\}\subseteq \cF\cup \{G\}$. If $3\le |G|\le n-1$, then either $\{1,2\}\subseteq G$ and there exists $3\le u\le n$ with $u \notin G$ or there exist $u\in \{1,2\}\setminus G$ and $\{v,w\}\subseteq G \cap [3,n]$. In the former case $2^{\{1,2,u\}}\setminus \{\{1,2\},\{1,2,u\}\}\cup \{G,[n]\}\subseteq \cF\cup \{G\}$ form a $\mathbb{B}_3$, while in the latter case  $2^{\{u,v,w\}}\setminus \{\{v,w\}.\{u,v,w\}\}\cup \{G,[n]\}\subseteq \cF\cup \{G\}$ form a $\mathbb{B}_3$.
\end{proof}

\textbf{Conflict of interest}

The authors declare that they have no conflict of competing interest.

\medskip

\textbf{Funding}

Shengjin  Ji is supported by  Natural Science Foundation of Shandong Province, China (No. ZR2022MA077) and  Postgraduate Education Reform Project of Shandong Province, China (No. SDYKC2023107).

Erfei Yue is supported by ERC Advanced grant GeoScape, No. 882971.

\end{document}